\documentclass[12pt,preprint,epsfig]{aastex}

\begin{document}
\title{Asymptotic harmonic behavior in the prime number distribution}
\author{Maurice H.P.M. van Putten}
\affil{Yongsil-Gwan 614, Sejong University, 143-747 Seoul, South Korea}

\begin{abstract}
We consider $\Phi(x)=x^{-\frac{1}{4}}\left[1-2\sqrt{x}\Sigma e^{-p^2\pi x}\ln p\right]$ on $x>0$, where the sum is over all primes $p$. 
If $\Phi$ is bounded on $x>0$, then the Riemann hypothesis is true or there are infinitely many zeros Re~$z_k>\frac{1}{2}$. 
The first 21 zeros give rise to asymptotic harmonic behavior in $\Phi(x)$ defined by the prime numbers up to one trillion. 
\end{abstract}

\section{Introduction}

The Riemann-zeta function is the analytic extension of
\begin{eqnarray}
\zeta(z)=1+\frac{1}{2^z}+\frac{1}{3^z} + \cdots = \Pi \left(1-p^{-z}\right)^{-1}~~(\mbox{Re}\,z>1),
\label{EQN_Ze}
\end{eqnarray}
where Euler's identity on the right hand side expresses the relation of the integers to the primes.
The zeros $z_k$ of Riemann's analytic continuation of (\ref{EQN_Ze}) comprise the negative even integers, $-2,-4,\cdots$, and an infinite number of nontrivial zeros $z_k=a_k+iy_k$ in the strip $0<a_k<1$. 

\begin{figure}[h]
\centerline{\includegraphics[scale=0.45]{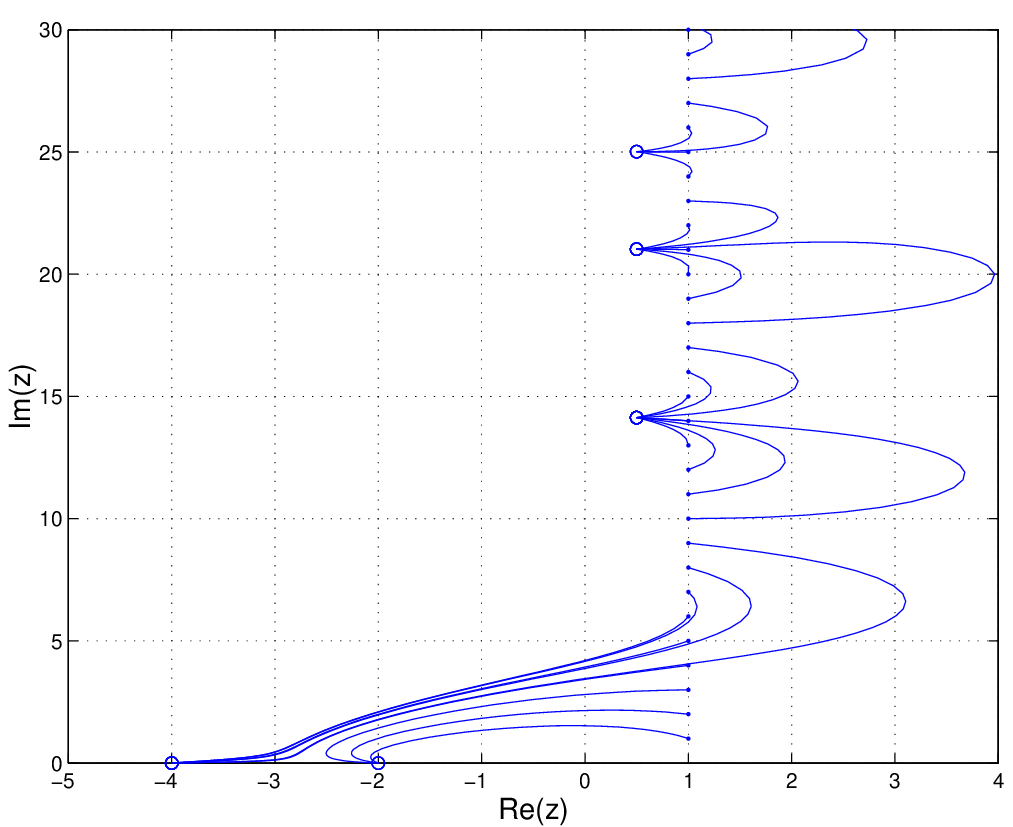}}
\caption{Shown are the trajectories of continuation $z(\lambda)$ in the complex plane $z$ by numerical integration of (\ref{EQN_C}) with initial data $z_0=1+ni$ $(n=1,2,3,\cdots)$ indicated by small dots on Re($z$)=1. Continuation produces roots indicated by open circles, 
defined by finite endpoints of $z(\lambda)$ in the limit as $\lambda$ approaches infinity. The roots produced by the choice of initial data are the first three on Re $z=\frac{1}{2}$ and -2 and -4 of the trivial roots.}
\label{fig_3}
\end{figure}

A general approach to find zeros is by continuation \citep{kel87}. If $z(0)=z_0$ is a starting point of a path $z(\lambda)$  with tangent $\tau=z^\prime(\lambda)$,
\begin{eqnarray}
\tau \frac{\zeta^\prime(z)}{\zeta(z)}=-1,
\label{EQN_C}
\end{eqnarray}
then the endpoint $z_*=\lim_{\lambda\rightarrow\infty} z(\lambda)$ is a zero of $\zeta(z)$, all of which are isolated. 
All known nontrivial zeros satisfy $\mbox{Re}~z_k=\frac{1}{2}$ to within numerical precision, the first three of which 
are $z_1=\frac{1}{2}\pm 14.1347i$, $z_2=\frac{1}{2}\pm 21.0220i,$ $z_3=\frac{1}{2}\pm 25.0109i.$ By the symmetry
\begin{eqnarray}
\zeta(s)=\zeta(1-s)\frac{\chi(1-s)}{\chi(s)},~~\chi(s) =\pi^{\frac{s}{2}} \Gamma\left(\frac{s}{2}\right),~~\Gamma(z) = \int_0^\infty t^{z-1} e^{-t} dt,
\label{EQN_ID}
\end{eqnarray}
it suffices to study zeros in the half plane Re($z)\ge\frac{1}{2}$. Fig. \ref{fig_3} illustrates root finding by (\ref{EQN_C}) for the first few zeros.

Continuation (\ref{EQN_C}) is determined by the prime numbers, since
\begin{eqnarray}
-\frac{\zeta^\prime(z)}{\zeta(z)} = - \sum \frac{\ln p}{p^z-1} = \Sigma \xi(mz),~~\xi(z)=\sum p^{-z} \ln p~~(\mbox{Re}\,z>1),
\end{eqnarray} 
whereby 
\begin{eqnarray}
\xi(z) = -\frac{\zeta^\prime(z)}{\zeta(z)}  - \sum_{m\ge 2} \xi(mz).
\label{EQN_B2}
\end{eqnarray}
The poles of $\xi(z)$ at the zeros are therefore expressed by the prime number distribution. 

In this paper, we study the distribution of zeros $z_k$ by Fourier analysis of the function
\begin{eqnarray}
\Phi(x)  = x^{-\frac{1}{4}}\left[ 1-2\sqrt{x}\phi(x)\right]
\label{EQN_E}
\end{eqnarray}
on $x>0$, where 
 \begin{eqnarray}
\varphi(x)=\sum  e^{-p^2\pi x}\log p
\label{EQN_D}
\end{eqnarray}
with summation over all primes. In what follows, we put
\begin{eqnarray}
Z(\lambda) = \sum \alpha_k e^{-\lambda(z_k-\frac{1}{2})},~~\alpha_k = \gamma(z_k),~~\gamma(z)=\frac{\Gamma\left(\frac{z}{2}\right)}{\pi^\frac{z}{2}}.
\label{EQN_Z0}
\end{eqnarray}
The $Z_k$ are absolutely summable by Stirling's formula and the asymptotic distribution of $z_k$.

{\bf Theorem 1.1.} {\em In the limit as $x>0$ becomes small, we have the asymptotic behavior}
\begin{eqnarray}
\Phi(x) =  \frac{1}{2} \gamma\left(\frac{1}{2}\right) +Z\left(\ln\sqrt{x}\right)+\frac{1}{3} \gamma\left(\frac{1}{3}\right) x^\frac{1}{12} 
+o\left( x^\frac{1}{12}\right).
\label{EQN_LA5}
\end{eqnarray}

In (\ref{EQN_LA5}), $Z$ is evidently unbounded in the limit as $x$ approaches zero whenever a finite number of zeros $z_k$ exists off the critical line Re\,$z=\frac{1}{2}$.

{\bf Corollary 1.2.} {\em If $\Phi$ is bounded, then the Riemann hypothesis is true or there are infinitely many zeros Re~$z_k>\frac{1}{2}$.} 

A similar relation between the distribution of $z_k$ and the primes is \citep{had93,von95}
\begin{eqnarray}
\frac{u-\psi_C(u)}{\sqrt{u}} = \sum \frac{u^{z_k-\frac{1}{2}}}{z_k} + \frac{\ln (2\pi)+ \ln\sqrt{1-u^{-2}}}{\sqrt{u}}
\label{EQN_vm}
\end{eqnarray}
based on the Chebyshev functions 
\begin{eqnarray}
\psi_C(u)=\sum_{p^k\le u}\ln (p),~~\vartheta_C(u)=\sum_{p\le u}\ln p,
\end{eqnarray}
where the sum is over all primes $p$ and integers $k$. In (\ref{EQN_LA5}), $\Phi(x)$ has a normalization by $x^\frac{1}{4}$ according to  and $Z$ is absolutely convergent for all $x>0$,  whereas in (\ref{EQN_vm}) $\psi_C(u)$ is normalized by $\sqrt{u}$ and the sum $\Sigma \frac{u^{z_k-\frac{1}{2}}}{z_k}$ is not absolutely convergent. Similar to Corollary 1.2, the left hand side of (\ref{EQN_vm}) will be bounded in the limit of large $u$ if the Riemann hypothesis is true.

\S2 presents some preliminaries on $\zeta(z)$. \S3 gives an integral representation of $\zeta(z)$ and a discussion on its singularity at $z=1$.  In \S4, Cauchy's integral formula is applied to derive a sum of residues associated with the $z_k$. The proof Theorem 1.1 follows from a Fourier transform and asymptotic analysis  (\S5). In \S6, we illustrate a direct evaluation of $\Phi(x)$ using the primes up to one trillion, showing harmonic behavior arising from $Z$ by the first few zeros $z_k$. We summarize our findings in \S7.

\section{Background}

Our analysis begins with some known properties of $\zeta(z)$ in, e.g., \cite{tit86,leh88,dus99,kei92,for02}). 

Riemann obtained an analytic extension of $\zeta(z)$ by expressing $n^{-z}$ in terms of $\Gamma\left(\frac{z}{2}\right)$,
\begin{eqnarray}
\gamma(z)\zeta(z) = \int_0^\infty x^{\frac{z}{2}-1} \theta_1(x)  dx,
\label{EQN_I1}
\end{eqnarray}
where
\begin{eqnarray}
\theta_1(x)= \frac{\theta(x)-1}{2},~~\theta(x)=\sum_{n=-\infty}^{\infty} e^{-n^2\pi x}. 
\label{EQN_TH}
\end{eqnarray}
Here, $\theta_1(x)$ satisfies $\theta_1(x)\sim \frac{1}{2\sqrt{x}}$ as $x$ approaches zero by the identity $\theta(x^{-1})=\sqrt{x}\theta(x)$ for the Jacobi function $\theta(x)$.
\footnote{When $z=n$ is an integer, $\frac{\pi^\frac{n}{2}}{\Gamma\left(\frac{n}{2}\right)}$ is one-half the surface area of $S^n$.} On Re$~z>1$, it obtains the meromorphic expression (e.g. Borwein et al. 2006)
\begin{eqnarray}
\gamma(z)\zeta(z) = \frac{1}{z(z-1)}+f(z),~~
f(z)=\int_1^\infty \left(x^{\frac{z}{2}-1}+x^{-\frac{z}{2}-\frac{1}{2}}\right) \theta_1(x) dx,
\label{EQN_R}
\end{eqnarray}
which gives a maximal analytic continuation of $\zeta(z)$ and shows a simple pole at $z=1$ with residue 1.

Riemann further introduced the symmetric form $Q(z)\zeta(z)$, 
$Q(z)=\frac{1}{2}z(z-1)\gamma(z)$ satisfying $Q(z)\zeta(z)=Q(1-z)\zeta(1-z)$, whereby
\begin{eqnarray}
\zeta(z)=\pi^{z-1}\zeta(1-z) \frac{\Gamma(\frac{1}{2}-\frac{z}{2})}{\Gamma(\frac{z}{2})} 
             =\frac{\pi^{z}\zeta(1-z) }{\cos\left(\frac{1}{2}\pi z\right)\Gamma(\frac{z}{2})\Gamma(\frac{1}{2}+\frac{z}{2})}
             =\frac{\pi^{z-\frac{1}{2}}{2^{z-1}}\zeta(1-z) }{\cos\left(\frac{1}{2}\pi z\right)\Gamma(z)}
\end{eqnarray}
using $\Gamma(\frac{1}{2}-\frac{z}{2})\Gamma(\frac{1}{2}+\frac{z}{2})=\frac{\pi}{\cos(\pi z)}$ and $\Gamma(z)\Gamma(z+\frac{1}{2})=2^{1-2z}\sqrt{\pi}\Gamma(2z)$. 
Along $z=1-iy$, $\zeta(z)$ is non-vanishing \citep{lit22,lit24,lit27,win41}, allowing 
\begin{eqnarray}
\frac{\zeta^\prime(z)}{\zeta(z)}=-\frac{\zeta^\prime(1-z)}{\zeta(1-z)}+\ln(2\pi) +\frac{\pi}{2}\tan\left(\frac{\pi z}{2}\right)-\psi(z)
\label{EQN_lk}
\end{eqnarray}
in terms of the digamma function
\begin{eqnarray}
\psi(z)=\frac{\Gamma^\prime(z)}{\Gamma(z)} \sim \ln(z) +O(z^{-1})
\label{EQN_psi}
\end{eqnarray}
in the limit of large $|z|$.

{\bf Lemma 2.1.} {\em In the limit of large $y$, the logarithmic derivative of $\zeta(z)$ satisfies}
\begin{eqnarray}
\frac{\zeta^\prime(iy)}{\zeta(iy)} = - \frac{\zeta^\prime(1-iy)}{\zeta(1-iy)} + O(\ln y ).
\end{eqnarray} 

{\em Proof.} The result follows from (\ref{EQN_psi}) and (\ref{EQN_lk}). $\Box$

{\bf Lemma 2.2.} {\em Along the line $z=iy$, we have the asymptotic expansion} $\left|\gamma(iy)\right| \sim \sqrt{\frac{2\pi}{y}} e^{-\frac{\pi}{2}y}$ {\em in the limit of large $y$, whereby the $\alpha_k$ are absolutely summable.}

{\em Proof.} Recall (\ref{EQN_Z0}) and the asymptotic expansion $\Gamma(z) = \sqrt{2\pi} z^{z-\frac{1}{2}} e^{-z}\left[ 1+O(z^{-1})\right]$ with a branch cut along the negative real axis. In the limit of large $y_k$, $y_k\sim \frac{2\pi k}{\ln k}$, and hence $\left|Z_k\right| \sim  e^{-\frac{\pi^2 k }{\ln k}}$, since $|\arg z_k|\sim \frac{\pi}{2}$ as $k$ becomes large. Hence, the $Z_k$ are absolutely summable. Numerically, their sum is small, $\Sigma|\alpha_k|=3.5\times 10^{-5}$ based on a large number of known zeros $z_k$. $\Box$

{\bf Lemma 2.3.} {\em In the limit of large $y$, we have}
\begin{eqnarray}
\left| \gamma(iy) \frac{\zeta^\prime(iy)}{\zeta(iy)}\right|  =  O\left( y^{-\frac{1}{2}}e^{-\frac{\pi}{2}y}\ln y \right). 
\end{eqnarray}
{\em Proof.} By Lemma 2.1-2, we have
\begin{eqnarray}
\left| \gamma(iy) \frac{\zeta^\prime(iy)}{\zeta(iy)}\right| 
\sim \sqrt{\frac{2\pi}{y}}\left( \left| \frac{\zeta^\prime(1-iy)}{\zeta(1-iy)}\right|  +   O(\ln y) \right) e^{-\frac{\pi}{2}y}
\end{eqnarray}
for large $y$. Also \citep{ric67,che99,tit86,bra05}
\begin{eqnarray}
\left| \frac{\zeta^\prime(1-iy)}{\zeta(1-iy)}\right| \le c  (\ln y)^\frac{2}{3} (\ln\ln y)^\frac{1}{3} 
\end{eqnarray}
on $y>\delta$ for some positive constants $c, \delta$. $\Box$

\section{An integral representation of $\xi(z)$}

Following the same steps leading to the Riemann integral for $\zeta(z)$, we have
\begin{eqnarray}
\gamma(z)\xi(z) = \int_0^\infty x^{\frac{z}{2}-1}\varphi(x)  dx = \frac{1}{z-1} + g(z),
\label{EQN_DX}
\end{eqnarray}
where $1/(z-1)$ absorbs the simple pole in $\xi(z)$ at $z=1$ due to the simple pole in $\zeta(z)$ at $z=1$,
leaving $g(z)$ analytic at $z=1$. Following a decomposition $g(z)=g_2(z)-g_1(z)$,  
\begin{eqnarray}
g_1(z)=\frac{1}{2}\int_0^1x^{\frac{2z-1}{4}}\Phi(x)\frac{dx}{x},~~ g_2(z)=\int_1^\infty x^{\frac{z}{2}-1}\varphi(x)dx,
\label{EQN_R4a}
\end{eqnarray}
and substitution $x=e^{2\lambda}$, $g(z)$ appears as the Laplace transforms
\begin{eqnarray}
g_1(z)=\int_{-\infty}^0 \Phi(e^{2\lambda})e^{\lambda (z-\frac{1}{2})} d\lambda,~~g_2(z)  = 2 \int_0^\infty \varphi(e^{2\lambda})e^{\lambda z} d\lambda.
\label{EQN_R4}
\end{eqnarray}
These integral expressions allow continuations to $\mbox{Re}~z>1$, respectively, the entire complex plane.

{\bf Lemma 3.1.} {\em Analytic extension of $g_1(z)$ extends to $z>\frac{1}{2}$.}

{\em Proof.} With $z=a+ib$, the second term on the right hand side in (\ref{EQN_B2}) satisfies
\begin{eqnarray}
\sum_{m\ge3} |\xi(mz)|\le\sum_{n\ge3} \frac{n^{-3a}\log n}{1-n^{-a}} < - \frac{\sqrt{2}\zeta^\prime(3a)}{\sqrt{2}-1},
\end{eqnarray} 
which is bounded in Re~$z=a>\frac{1}{2}$. 
Since the second term $\xi(2z)$ in (\ref{EQN_B2}) is analytic in Re~$z=a>\frac{1}{2}$, it follows that $g(a)$ in is analytic on $a>\frac{1}{2}$. Following (\ref{EQN_B2}) as $a$ approaches $\frac{1}{2}$ from the right, we have
\begin{eqnarray}
\xi(a) = -\frac{1}{2a-1}+u_1(a),
\label{EQN_A1}
\end{eqnarray}
where $u_1(a)$ is analytic at $a=\frac{1}{2}$. By (\ref{EQN_DX}), as $a$ approaches $\frac{1}{2}$ from the right, we have 
\begin{eqnarray}
g_1(a)= - \frac{1}{2a-1} +u_2(a),
\label{EQN_A2}
\end{eqnarray}
where $u_2(a)$ is analytic about $a=\frac{1}{2}$. $\Box$

\begin{figure}[h]
\centerline{\includegraphics[scale=0.4]{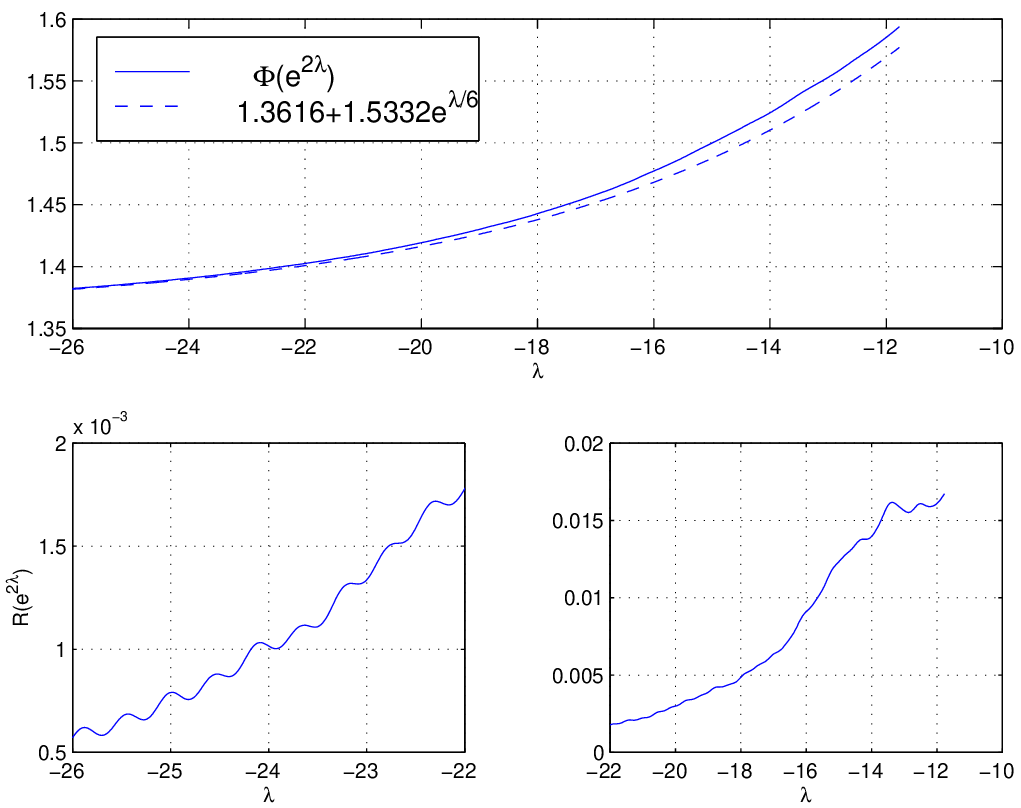}}
\caption{
The top window shows ${\Phi}\left(e^{2\lambda}\right)$ on $\lambda\epsilon[-26, -11.7756]$ and its leading order approximation $1.3616+1.5332e^{\frac{\lambda}{6}}$. The asymptotic harmonic behavior is apparent in the residual difference (\ref{EQN_Rx}) between the two, shown in the bottom two windows, including the period of $2.2496$ in $\lambda$ associated with the first zero $z_*=\frac{1}{2}\pm 14.1347i$.}
\label{fig_3d}
\end{figure}

Fig. \ref{fig_3d} shows a numerical evaluation of ${\Phi}(x)$ for small $x$ evaluated for the 37.6 billion primes up to one trillion,
allowing $x$ down to $2.6\times 10^{-23}$ ($\lambda=-26$) in view of the requirement for an accurate truncation in $\varphi(x)$ as defined 
by (\ref{EQN_D}). The result shows asymptotic harmonic behavior in the limit as $x$ becomes small.

If the integral 
\begin{eqnarray}
\int_\epsilon^1x^{\frac{2z-1}{4}}\Phi(x)\frac{dx}{x}
\label{EQN_R4c}
\end{eqnarray}
is absolutely convergent as $\epsilon>0$ approaches zero, e.g., when $\Phi(x)$ is of one sign in some neighborhood of $z=0$, as in
the numerical evaluation shown in Fig. \ref{fig_3d}, then $g_1(z)$ has an analytic extension into Re~$z>\frac{1}{2}$ with no singularities, implying the absence of $z_k$ in this region. 
However, this requires information on the point wise behavior of $\Phi(x)$, which goes beyond the relatively weaker integrability property (\ref{EQN_R4a}). 

To make a step in this direction, we next apply a linear transform to (\ref{EQN_B2}) to derive the asymptotic behavior of $\Phi(x)$ in terms of the distribution $z_k$.

\section{A sum of residues $Z$ associated with the non-trivial zeros}

Consider
\begin{eqnarray}
h(z) = \gamma(z) \frac{\zeta^\prime(z)}{\zeta(z)}+\frac{1}{z-1}
\label{EQN_gk}
\end{eqnarray}
and its Fourier transform
\begin{eqnarray}
H(\lambda)=\int_{a-i\infty}^{a+i\infty} h(z)e^{-\lambda z} \frac{dz}{2\pi i}.
\label{EQN_EH1}
\end{eqnarray}

{\bf Lemma 4.1.} {\em $h(z)$ has a simple pole at $z=1$ with residue 1 and simple poles at each of the nontrivial zeros $z_k$ of $\zeta(z)$ with residue $Z_k$.}

{\em Proof.}  We have (e.g. Borwein et al. 2006)
\begin{eqnarray}
\zeta_1(z)=\frac{1}{2}z(z-1)\gamma(z)\zeta(z),
~~\frac{\zeta_1^\prime(z)}{\zeta_1(z)}=B + \Sigma_k \left( \frac{1}{z-z_k}+\frac{1}{z_k}\right),
\label{EQN_zet1}
\end{eqnarray}
where $B$ is a constant, so that
\begin{eqnarray}
\gamma(z)\frac{\zeta^\prime(z)}{\zeta(z)}+\frac{1}{z-1}=
    \gamma(z)\left[B + \Sigma_k \left( \frac{1}{z-z_k}+\frac{1}{z_k}\right)\right]+A(z).
\label{EQN_h}
\end{eqnarray}
Here
\begin{eqnarray}
A(z)=\frac{1-\gamma(z)}{z-1}-\frac{\gamma(z)}{z} - 2\frac{\gamma(z)\psi(z)}{\ln \pi},
\end{eqnarray}
where $\psi(z)$ denotes the digamma function as before, 
includes contributions from the logarithmic derivative of the factor to $\zeta(z)$ in (\ref{EQN_zet1}), whose singularities are
restricted to the trivial zeros of $\zeta(z)$. $\Box$

We now consider the Fourier integral over Re\,$z=a$ as part of contour integration closed over $z=x\pm iY$ and Re\,$z=0$.

{\bf Proposition 4.2.} {\em The Fourier transform of $h(z)$ over Re~$z>\sup a_k$ satisfies}
\begin{eqnarray}
H(\lambda)=e^{-\frac{\lambda}{2}}Z(\lambda)+O(1)
\label{EQN_EH2}
\end{eqnarray}
{\em in the limit of large $\lambda<0$.}

{\em Proof.} Integration over $z=x+iY$ $(0<x<a)$ gives
\begin{eqnarray}
e^{-iY}\int_{iY}^{iY+a} h(z)e^{-\lambda x}\frac{dx}{2\pi i} =
e^{-iY}\int_{iY}^{iY+a} \frac{\Gamma\left(\frac{z}{2}\right)}{\pi^{\frac{z}{2}}}\frac{\zeta^\prime(z)}{\zeta(z)}\frac{dx}{2\pi i}+O\left(Y^{-1}\right),
\label{EQN_Y}
\end{eqnarray}
where we choose $Y$ to be between two consecutive values of $y_k$. We have
\begin{eqnarray}
\int_{iY}^{iY+a}\frac{\Gamma\left(\frac{z}{2}\right)}{\pi^{\frac{z}{2}}}\frac{\zeta^\prime(z)}{\zeta(z)}\frac{dx}{2\pi i} \sim
{\gamma_k} \int_{iY}^{iY+a}\frac{\zeta^\prime(z)}{\zeta(z)}\frac{dx}{2\pi i}\sim \frac{\gamma_k}{2\pi i} \ln(2a-1) \left[1+\frac{4ai(y_k-Y)}{1-2a}+ \pi i\right].
\end{eqnarray}
In the limit as $k$ approaches infinity, $y_k-Y$ approaches zero and $|\gamma_k|$ becomes small by Lemma 2.2., whence
\begin{eqnarray}
\left(\int_{iY}^{iY+a}-\int_{-iY}^{-iY+a}\right) \frac{\Gamma\left(\frac{z}{2}\right)}{\pi^{\frac{z}{2}}}\frac{\zeta^\prime(z)}{\zeta(z)}\frac{dx}{2\pi i} \sim
 \ln(2a-1){\rm Im}~\gamma_k = O\left( \ln(2a-1) \sqrt{\frac{2\pi}{y_k}}e^{-\frac{\pi}{2}y_k}\right).
\end{eqnarray}
Next, integration over $z=iy$ with a small semicircle around $z=0$ obtains an $O(1)$ result in the limit of large $\lambda$
by application of Lemma 2.1-3 and the Riemann-Lebesque Lemma.
The result now follows in the limit as $k$ approaches infinity, taking into account the residue sum $e^{-\frac{\lambda}{2}}Z(\lambda)$ associated with the $z_k$ and absolute summability of the $\alpha_k$. $\Box$

\section{Proof of Theorem 1.1}

Multiplying (\ref{EQN_B2}) by $\gamma(z)$, we have
\begin{eqnarray}
 \gamma(z)\xi(z) = - \gamma(z) \frac{\zeta^\prime(z)}{\zeta(z)} - \gamma(z)\sum_{m\ge2}\xi(mz),
\end{eqnarray}
that is, by (\ref{EQN_DX}) and (\ref{EQN_gk}),
\begin{eqnarray}
\frac{1}{z-1} + g(z) =  - h(z) + \frac{1}{z-1} - \gamma(z)\sum_{m\ge2}\xi(mz).
\end{eqnarray}
We thus consider
\begin{eqnarray}
g_1(z)= g_2(z) + h(z) + \gamma(z)\sum_{m\ge2}\xi(mz),
\label{EQN_g1}
\end{eqnarray}
which {\em ab initio} is defined on Re\,$z>1$ by Euler's identity with Fourier transform
\begin{eqnarray}
G_1(\lambda)=\int_{a-i\infty}^{a+i\infty}g_1(z)e^{-\lambda z}\frac{dz}{2\pi i}  = e^{-\frac{\lambda}{2}}\Phi(e^{2\lambda}). 
\label{EQN_G1}
\end{eqnarray}
Turning to the right hand side of (\ref{EQN_g1}), we consider the coefficients
\begin{eqnarray}
c_m(z)=\frac{\gamma(z)}{\gamma(mz)},~~C_m=\frac{1}{m}\gamma\left(\frac{1}{m}\right) ~~(m\ge1). 
\label{EQN_cm}
\end{eqnarray}
Here, $C_m=m^{-1}\gamma(1/m)$ since $\gamma(1)=1$. In particular, $C_2=\frac{1}{2}\gamma\left(\frac{1}{2}\right)$ and $m^{-1}c_m(m^{-1}z)=1+\left(\frac{1}{2}\ln\pi+\frac{1}{2}\gamma\right)z+O\left(z^2\right)$ has a well defined limit and $C_m\rightarrow2$ in the limit as $m$ becomes arbitrarily large.

{\bf Lemma 5.1.} {\em The sum $\sum_{m\ge n} \xi(mz)$ is well-defined on Re~$z> \frac{1}{n}$.}

{\em Proof.} The result follows from the case $n=2$. By the Prime Number Theorem, $p_k\sim k\ln k$, whereby summation over the tails $k\ge n$ satisfy 
\begin{eqnarray}
\sum_{k\ge n}^\infty \frac{\ln(p_k)}{p_k^{2a}}\sim \sum_{k\ge n}^{\infty} \left[ \frac{1}{k^{2a}\ln(k)^{2a-1}}+\frac{\ln\ln(k)}{(k\ln(k))^{2a}}\right]<\infty
\end{eqnarray}
whenever $a>\frac{1}{2}$. Hence, for $z=a+iy$, $\left|\Sigma p^{-2z} \ln p \right|\le \Sigma p^{-2a} \ln p <\infty$ whenever $a>\frac{1}{2}$. It follows that 
\begin{eqnarray}
\left|\Sigma_{m\ge2}\xi(mz)\right|\le \Sigma_{m\ge2} \Sigma_p  p^{-(m-2)a} p^{-2a} \ln p \le \Sigma_{m\ge0} 2^{-m}\Sigma_p p^{-2a} \ln p =\Sigma_p p^{-2a} \ln p <\infty
\end{eqnarray}
on Re~$z>\frac{1}{2}$. $\Box$ 

{\bf Lemma 5.2.} {\em For any $m\ge 2$, the Fourier transform of $\frac{c_m(z)}{mz-1}$ over Re~$z=a>\frac{1}{2}$ satisfies}
\begin{eqnarray}
{D}_m(\lambda)=C_m e^{-\frac{\lambda}{m}} + o(1)
\label{EQN_EH2}
\end{eqnarray}

{\em Proof.} The Fourier integral can be obtained in a contour integration with closure over $z=iy$ and the edges $z=x+iY$ $(0<x<a)$ for large $\pm Y$. In the notation (\ref{EQN_cm}), it obtains a residue $C_m=m^{-1}c_m(1/m)=m^{-1}\gamma(1/m)$ at $z=1/m$, since $\gamma(1)=1$, whence
\begin{eqnarray}
{D}_m(\lambda)=C_m e^{-\frac{\lambda}{m}} +e^\frac{\lambda}{2} \int_{-\infty}^{\infty} \frac{c_m(iy)}{imy-1}e^{-i\lambda y} \frac{dy}{2\pi}.
\label{EQN_EH2}
\end{eqnarray}
The integral (\ref{EQN_EH2}) exists by virtue of a removable singularity of $c_m(z)$ at $z=0$. It asymptotically decays to zero for large $\lambda$ when $m\ge2$ by the Riemann-Lebesque Lemma. $\Box$

We now consider (\ref{EQN_g1}) with (\ref{EQN_DX}),
\begin{eqnarray}
g_1(z)= g_2(z) + h(z) + \sum_{m\ge2}\left( \frac{c_m(z)}{mz-1} + c_m(z)g(mz)\right)=
h(z)+ \sum_{m=2}^{N} \frac{c_m(z)}{mz-1}+r_N(z)
\label{EQN_fN}
\end{eqnarray}
with a remainder
\begin{eqnarray}
r_N(z)=g_2(z)+\sum_{m\ge2} c_m(z)g(mz)+ \gamma(z)\sum_{m\ge N+1}\xi(mz).
\label{EQN_RN}
\end{eqnarray}

{\bf Lemma 5.3.} {\em For $N\ge 3$, the Fourier transform}
\begin{eqnarray}
e^{\frac{\lambda}{2}}R_N(\lambda) =  \int_{a-i\infty}^{a+i\infty}  r_N(z) e^{-\lambda (z-\frac{1}{2})}\frac{dz}{2\pi i} = o(1) 
\label{EQN_RNL}
\end{eqnarray}
{\em in the limit of large $\lambda<0$.}  

{\em Proof.} Since $r_N(z)$ is analytic in Re~$z > \frac{1}{3}$, we are at liberty to consider the transform $e^{\frac{\lambda}{2}}R_N(\lambda)$ on $a=1/2$. The result follows from the Riemann-Lebesque Lemma. $\Box$

{\em Proof of Theorem 1.1.} The Fourier transform of (\ref{EQN_fN}) is
\begin{eqnarray}
G_1(\lambda)=H(\lambda)+D_2(\lambda)+D_3(\lambda)+R_3(\lambda).
\label{EQN_g4}
\end{eqnarray}
By Proposition 4.2 and Lemmas 5.1-5.2, we have 
\begin{eqnarray}
e^{-\frac{\lambda}{2}}\Phi(e^{2\lambda}) = e^{-\frac{\lambda}{2}}Z(\lambda) +C_2e^{-\frac{\lambda}{2}}+C_3e^{-\frac{\lambda}{3}}+o\left(e^{-\frac{\lambda}{2}}\right).
\label{EQN_g4b}
\end{eqnarray}
With $x=e^{2\lambda}$, Theorem 1.1 now follows. $\Box$

\section{Numerical illustration of asymptotic harmonic behavior}

\begin{figure}[h]
\centerline{\includegraphics[scale=0.6]{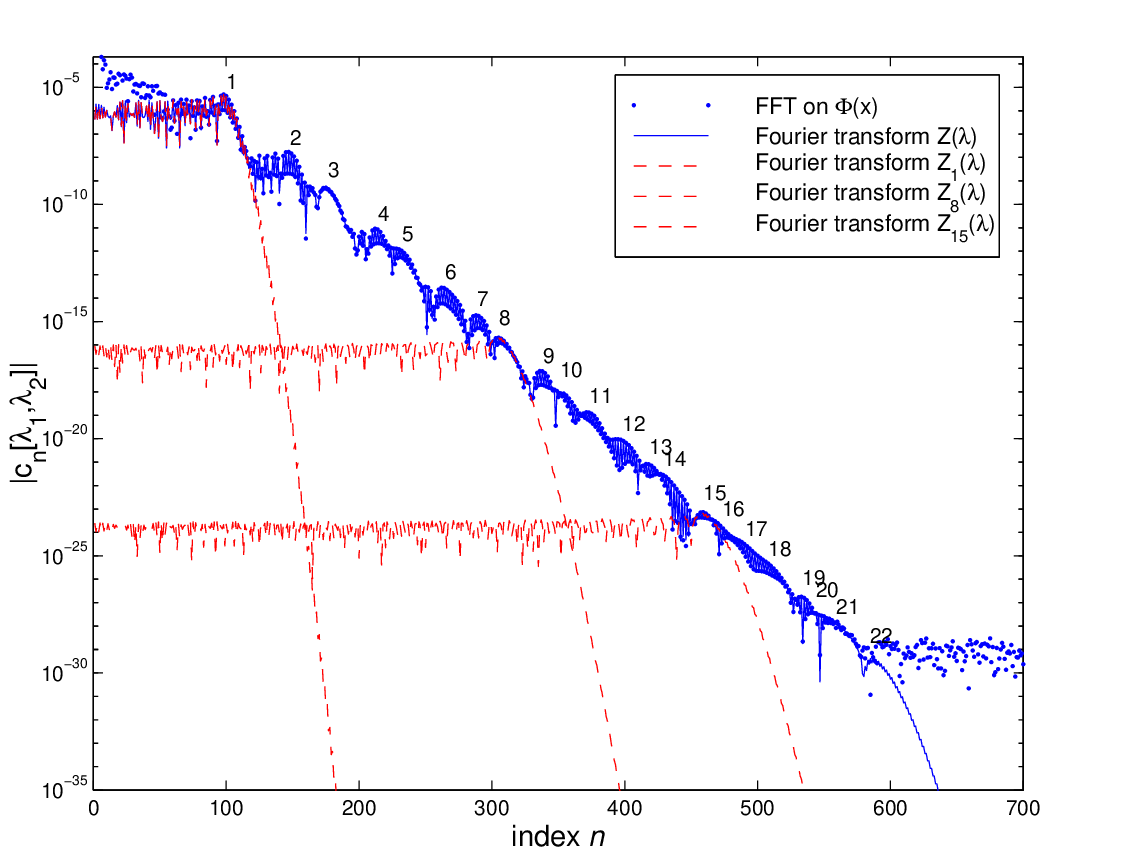}}
\caption{
Shown are the absolute values of the Fourier coefficients $c_n[\lambda_1,\lambda_2]$ of ${\Phi}\left(e^{2\lambda}\right)$ obtained by a Fast Fourier Transform (FFT) of (\ref{EQN_Rx}) on the computational domain (\ref{EQN_LAP}), where $\lambda_1=-26$, $\lambda_2= -11.7756]$ covers 32 periods of $Z_1(\lambda)$ ($dots$), on the basis of the 37,607,912,2019 primes up to 1,000,000,000,0039. The resulting spectrum is compared with the exact spectra $c_{ni}[\lambda_1,\lambda_2]$ of the $Z_i(\lambda)$ given by the analytic expression (\ref{EQN_cn}) for $i=1,2,3,\cdots$ ($continuous$ $line$). Shown are also the individual spectra of $Z_i(\lambda)$ for $i=1,8$ and 15 associated with the zeros $z_1$, $z_8$ and $z_{15}$. The match between the computed and exact spectra accurately identifies the first 21 harmonics of $Z(\lambda)$ in ${\Phi}$ out of 22 shown, corresponding to the first 21 nontrivial zeros $z_i$ of $\zeta(z)$.}
\label{fig_3e}
\end{figure}

The harmonic behavior emerges in 
\begin{eqnarray}
R(x)=-{\Phi}(x) - C_2 - C_3x^\frac{1}{12}.
\label{EQN_Rx}
\end{eqnarray}
To search for higher harmonics $Z_i(\lambda)$ associated with the zeros $z_i$ in $\lambda \epsilon [-26, -11.7759]=[\lambda_1-\lambda_2,\lambda_1+\lambda_2]$, we compare the spectrum of ${\Phi}\left(e^{2\lambda}\right)$ by taking a Fast Fourier Transform with respect to $\alpha$, 
\begin{eqnarray}
\lambda(\alpha) = \lambda_1+\lambda_2\cos\alpha~~(\alpha\epsilon [0,2\pi]),
\label{EQN_LAP}
\end{eqnarray}
and compare the results with an analytic expression for the Fourier coefficients of the $Z_i(\lambda)$ $(i=1,2,\cdots)$, 
\begin{eqnarray}
c_{ni}[\lambda_1,\lambda_2]=2\mbox{Re}\{(-i)^n\gamma_i e^{-i\lambda_1z_i}J_n\left(-\lambda_2z_i\right)\},
\label{EQN_cn}
\end{eqnarray}
where $J_n(z)$ denotes the Bessel function of the first of order $n$. Fig. \ref{fig_3e} shows the first 21 harmonics in our evaluation of ${\Phi}(x)$, which is about the maximum that can be calculated by direct summation in quad precision.

\section{Conclusions}

The zeros $z_k=a_k+iy_k$ of the Riemann-zeta function are endpoints of continuation, defined by an
expressed by a regularized sum $\Phi(x)$ over the prime numbers defined by (\ref{EQN_E}). 

The zeros $z_k$ of $\zeta(z)$ introduce asymptotic harmonic behavior in ${\Phi}\left(e^{2\lambda}\right)$ as a function of $\lambda<0$, defined by the sum $Z(\lambda)$ of residues of the $z_k$, shown in Figs. 2-3. Primes up to 4 billion are needed to identify the first 4 harmonics, up to 70 billion for the 10 and up to 1 trillion for the first 21. It appears that the prime number range scales approximately exponentially with the number of harmonics it contains. 

Theorem 1.1 describes a correlation between the distribution of the primes and the distribution of the nontrivial zeros $z_k$. Suppose there are a finite number of zeros $z_k$ in Re~$z>\frac{1}{2}$. We may then consider $k^*$ for which $a_{k^*}=\max a_k$ gives rise to dominant exponential growth in $Z(\lambda)$ in the limit as $\lambda<0$ becomes large. This observation leads to Corollary 1.2. $Z$ can remain bounded in $x>0$ only if the Riemann hypothesis is true, or if $Z(\lambda)$ remains fortuitously bounded as an infinite sum over $a_k>\frac{1}{2}$ with no maximum in $a<1$. Conversely, Riemann hypothesis implies
\begin{eqnarray}
\lim_{x\rightarrow0^+}\Phi(x)=\frac{1}{2}\gamma\left(\frac{1}{2}\right) \simeq 1.3616.
\end{eqnarray}

According to (\ref{EQN_LA5}) and our numerical calculation shown in Fig. 3, the number of primes relevant to the observed asymptotic harmonic behavior scales approximately exponentially with the number of zeros $z_k$. The zeros $z_k$ explored to large values by existing
numerical experiments hereby constrain the distribution of an exponentially large number of primes.

{\bf Acknowledgment.} The author gratefully acknowledges stimulating discussions with Fabian Ziltener and Anton F.P. van Putten.
Some of the manuscript was prepared at the Korea Institute for Advanced Study, Dongdaemun-Gu, Seoul. This research was supported in part by the National Science Foundation through TeraGrid resources provided by Purdue University under grant number TG-DMS100033. We specifically acknowledge the assistance of Vicki Halberstadt, Rich Raymond and Kimberly Dillman. The computations have been carried out using Lahey Fortran 95.

\end{document}